%% file: Modular_Poly_Final_LMSJCM.tex
\documentclass[11pt]{amsart}

\input{MTemplate}

\title{Computing Modular Polynomials}

\author{Denis Charles}
\thanks{2000 Mathematics Subject Classification. Primary: 11G18, 11Y16. Secondary: 14H52, 14G50}
\address{Department of Computer Science, University of Wisconsin-Madison, Madison, WI - 53706.}
\email{cdx@cs.wisc.edu}
\author{Kristin Lauter}
\address{Microsoft Research, One Microsoft Way, Redmond, WA - 98052.}
\email{klauter@microsoft.com}

\begin{document}
\maketitle
\section{Introduction}
The $\ell^{th}$ modular polynomial, $\phi_{\ell}(x,y)$, parameterizes pairs of elliptic curves with
a cyclic isogeny of degree $\ell$ between them. Modular polynomials provide the defining equations for modular curves, 
and are useful in many different aspects of computational number theory and cryptography.  For example, computations with
modular polynomials have been used to speed elliptic curve point-counting algorithms~(\cite{bss99} Chapter VII).  \\


The standard method for computing modular polynomials consists of computing the Fourier expansion of 
the modular $\mathbf{j}$-function and solving a linear system of equations
to obtain the integral coefficients of $\phi_{\ell}(x,y)$. 
According to Elkies (see \cite{elk98} \S3) this method has a running time
of $O(\ell^{4+\epsilon})$ if one uses fast multiplication. However, our analysis (given in the Appendix) shows that
the running time of this method is, in fact, $\Theta(\ell^{9/2 + \epsilon})$ using fast multiplication. 
\\


The object of the current paper is to compute the modular polynomial (for prime $\ell$) directly modulo a prime $p$, without first computing the
coefficients as integers.  Once the modular polynomial has been computed for enough small primes, our approach can also 
be combined with the Chinese Remainder Theorem (CRT) approach as in~ \cite{cnst98} or \cite{alv03} to obtain the modular 
polynomial with integral coefficients or with coefficients modulo a much larger prime using Explicit CRT. 
Our algorithm does not involve computing Fourier coefficients of modular functions. The running time of our algorithm
turns out to be $O(\ell^{4+\epsilon})$ using fast multiplication. We believe our method is interesting as it is asymptotically faster; 
and is an essentially different
approach to computing modular polynomials. Furthermore, our algorithm also yields as a corollary a fast way to compute a
random $\ell$-isogeny of an elliptic curve over a finite field.\\

The idea of our algorithm is as follows.  Mestre's algorithm, M{\'e}thode des graphes~\cite{mes86}, uses the $\ell^{th}$ modular polynomial
modulo $p$ to navigate around the connected graph of supersingular elliptic curves over $\mathbb{F}_{p^2}$ in order to 
compute the number of edges (isogenies of degree $\ell$) between each node.  From the graph, Mestre then obtains the 
$\ell^{th}$ Brandt matrix giving the action of the $\ell^{th}$ Hecke operator on modular forms of weight $2$.  In our 
algorithm we do the opposite: we compute the $\ell^{th}$ modular polynomial modulo $p$ by computing all the isogenies 
of degree $\ell$ between supersingular curves modulo $p$ via V\'elu's formulae.  Specifically, for a given $\mathbf{j}$-invariant, $j$ (say), 
of
a supersingular elliptic curve over $\mathbb{F}_{p^2}$, Algorithm 1 computes $\phi_{\ell}(x,j)$ modulo $p$ by computing 
the $\ell +1$ distinct  subgroups of order $\ell$ and computing the $\mathbf{j}$-invariants of the $\ell +1$ corresponding 
$\ell$-isogenous elliptic curves. Algorithm 2 then uses the connectedness of the graph of supersingular elliptic curves 
over $\mathbb{F}_{p^2}$ to move around the graph, calling Algorithm 1 for different values of $j$ until enough information
is obtained to compute $\phi_{\ell}(x,y)$ modulo $p$ via interpolation.\\


There are several interesting aspects to Algorithms 1 and 2.  Algorithm 1 does not use the factorization of 
the $\ell$-division polynomials to produce the subgroups of order $\ell$.  Instead we generate independent $\ell$-torsion
points by picking random points with coordinates in a suitable extension of $\mathbb{F}_p$ and taking a scalar multiple
which is the group order divided by $\ell$.  This turns out to be more efficient than factoring the $\ell^{th}$
division polynomial for large $\ell$.  \\

Algorithm 2 computes $\phi_{\ell}(x,y)$ modulo $p$ by doing only computations with supersingular elliptic curves in
characteristic $p$ even though $\phi_{\ell}(x,y)$ is a general object giving information about isogenies between 
elliptic curves in characteristic $0$ and ordinary elliptic curves in characteristic $p$.  The advantage that we gain by 
using supersingular elliptic curves is that we can show that the full $\ell$-torsion is 
defined over an extension of degree $O(\ell)$ of the base 
field $\mathbb{F}_{p^2}$, whereas in general the field of definition can be of degree as high as $\ell^2 -1$.\\

In this article we provide a running time analysis assuming fast multiplication implementation of field operations.
But for small values of $\ell$ fast multiplication is usually not used in practice, thus we also give the running time (without the
analysis) assuming a na{\"i}ve implementation of field operations. 

\section{Local computation of $\phi_{\ell}(x,j)$}
The key ingredient of the algorithm is the computation of the univariate polynomial $\phi_{\ell}(x,j)$ modulo a prime $p$ given
a $\mathbf{j}$-invariant $j$. We describe the method to do this here. \\

{\bf Algorithm 1}\\
{\bf Input: } Two distinct primes $p$ and $\ell$, and $j$ the $\mathbf{j}$-invariant of a supersingular elliptic curve $E$
	over a finite field $\mathbb{F}_q$ of degree at most $2$ over a prime field of characteristic $p$.\\
{\bf Output: } The polynomial $\phi_{\ell}(x,j) = \prod_{E' \text{ $\ell$-isogenous to }E}(x - j(E')) \in \mathbb{F}_{p^2}[x]$.\\

\begin{enumerate}
\item[Step 1] Find the generators $P$ and $Q$ of $E[\ell]$:
	\begin{enumerate}
	\item Let $n$ be such that $\mathbb{F}_q(E[\ell]) \subseteq \mathbb{F}_{q^n}$.
	\item Let $S = \sharp E(\mathbb{F}_{q^n})$, the number of $\mathbb{F}_{q^n}$ rational points on $E$.
	\item Set $s = S/\ell^k$, where $\ell^k$ is the largest power of $\ell$ that divides $S$ (note $k \geq 2$).
	\item Pick two points $P$ and $Q$ at random from $E[\ell]$:
		\begin{enumerate}
		\item Pick two points $U,V$ at random from $E(\mathbb{F}_{q^n})$.
		\item Set $P' = s U$ and $Q' = s V$, if either $P'$ or $Q'$ equals $\mathcal{O}$ then repeat step (i).
		\item Find the smallest $i_1, i_2$ such that $\ell^{i_1} P' \neq \mathcal{O}$ and $\ell^{i_2} Q' \neq \mathcal{O}$
			but $\ell^{i_1 + 1} P' = \mathcal{O}$ and $\ell^{i_2 + 1}Q' = \mathcal{O}$.
		\item Set $P = \ell^{i_1} P'$ and $Q = \ell^{i_2} Q'$.
		\end{enumerate}
	\item Using Shanks's Baby-steps-Giant-steps algorithm check if $Q$ belongs to the group generated by $P$. If so
		repeat step (d).
	\end{enumerate}
\item[Step 2] Find the $\mathbf{j}$-invariants $j_1, \cdots, j_{\ell+1}$ in $\mathbb{F}_{p^2}$ 
of the $\ell+1$ elliptic curves that are $\ell$-isogenous to $E$.
	\begin{enumerate}
	\item Let $G_1 = \langle Q \rangle$ and $G_{1+i} = \langle P + (i-1)Q \rangle$ for $1 \leq i \leq \ell$.
	\item For each $i$, $1 \leq i \leq \ell+1$ compute the $\mathbf{j}$-invariant of the elliptic curve $E/G_i$
		  using V{\'e}lu's formulas.
	\end{enumerate}
\item[Step 3] Output $\phi_{\ell}(x,j) = \prod_{1 \leq i \leq \ell+1}(x - j_i)$.\\
\end{enumerate}

\begin{Rem} In step (1e), one could alternately use the Weil pairing to check whether $P$ and $Q$ 
generate the $\ell$-torsion. Doing so, however, does not lead to an asymptotic improvement in the
running time of the algorithm.
\end{Rem}

The following lemma gives the possibilities for the value of $n$ in Step (1a). We prove the following
result for all elliptic curves not just supersingular ones.

\begin{Lem} \label{lem_deg_bd:this}
Let $E/\mathbb{F}_q$ be an elliptic curve, and let $\ell$ be a prime not equal to the characteristic
of $\mathbb{F}_q$. Then $E[\ell] \subseteq E(\mathbb{F}_{q^n})$
where $n$ is a divisor of either $\ell(\ell-1)$ or $\ell^2 - 1$.
\end{Lem}
\begin{Proof}
The Weil-pairing tells us that if $E[\ell] \subseteq \mathbb{F}_{q^n}$ then $\mathbb{\mu}_{\ell} \subseteq \mathbb{F}_{q^n}$ (\cite{sil86} Corollary 8.1.1).
We seek, however, an upper bound on $n$, to do this we use the Galois representation coming from the $\ell$-division points of $E$.
Indeed, we have an {\em injective} group homomorphism (\cite{sil86} Chapter III, \S7)
\begin{align*}
	\rho_{\ell} : \Gal(\mathbb{F}_{q}(E[\ell])/\mathbb{F}_q) \rightarrow \mathrm{Aut}(E[\ell]) \cong \GL_2(\mathbb{F}_{\ell}).
\end{align*}
The Galois group $\Gal(\mathbb{F}_q(E[\ell])/\mathbb{F}_q)$ is cyclic, thus by $\rho_{\ell}$ the possibilities for 
$\Gal(\mathbb{F}_q(E[\ell])/\mathbb{F}_q)$ are limited
to cyclic subgroups of $\GL_2(\mathbb{F}_{\ell})$. In other words, we are interested in the orders
of the elements in $\GL_2(\mathbb{F}_{\ell})$. The elements of $\GL_2(\mathbb{F}_{\ell})$
are conjugate to one of the following types of matrices:
\begin{align*}
\begin{pmatrix}
\alpha & 0\\
0 & \beta
\end{pmatrix},
\begin{pmatrix}
\alpha & 1 \\
0 & \alpha
\end{pmatrix}, \text{ for } \alpha, \beta \in \mathbb{F}_{\ell}^*, \alpha \neq \beta,
\end{align*}
or those corresponding to multiplication by an element of $\mathbb{F}_{\ell^2}^*$
on the $2$-dimensional $\mathbb{F}_{\ell}$ vector space $\mathbb{F}_{\ell^2}$.
It is easy to see that the orders of these elements all divide $\ell(\ell-1)$ or $\ell^2-1$.
Thus the degree of the field extension containing the $\ell$-torsion points on $E$ must divide either $\ell(\ell-1)$ or $\ell^2-1$.
\end{Proof}

We will try step (1) with $n = \ell^2 - 1$, if steps (1d - 1e) do not succeed for some $K$ (a constant) many trials, we repeat
with $n = \ell(\ell - 1)$. The analysis that follows shows that a sufficiently large constant $K$ will work.\\

For step (1b) we do not need a point counting algorithm to determine $S$. Since $E$ is a supersingular elliptic curve, 
we have the following choices for the trace of Frobenius $a_q$:
\begin{align*}
	a_q &= \begin{cases}
		0 & \text{ if $E$ is over } \mathbb{F}_p \\
		0, \pm p, \pm 2p &\text{ if $E$ is over } \mathbb{F}_{p^2}.
		\end{cases}
\end{align*}
Not all the possibilities can occur for certain primes, but we will not use this fact here (see \cite{sch87}). 
If the curve is over $\mathbb{F}_{p^2}$ we can determine probabilistically the value of $a_q$ as follows. Pick a point
$P$ at random from $E(\mathbb{F}_q)$ and check if $(q + 1 + a_q)P = \mathcal{O}$. Since the pairwise gcd's of the possible
group orders divide $4$, with high probability only the correct
value of $a_q$ will annihilate the point. Thus in $O(\log^{2+o(1)} q )$ time we can determine with high probability
the correct value of $a_q$. Once we know the correct trace $a_q$, we can find the roots (in $\overline{\mathbb{Q}}$), 
$\pi$ and $\overline{\pi}$, of
the characteristic polynomial of the Frobenius $\phi^2 - a_q \phi + q$. Then the number of points lying on $E$ over 
the field $\mathbb{F}_{q^n}$ is given by $q^n + 1 - \pi^n - \overline{\pi}^n$, this gives us $S$.\\

{\bf Note: } We could have used a deterministic point counting algorithm to find $\sharp E(\mathbb{F}_{q})$ but this
would have cost $O(\log^6 q)$ field operations.\\

A lower bound on the probability that step (1d) succeeds is given by the following lemma whose proof
is straightforward.

\begin{Lem}\label{lem_prb_suc:this}
For a random choice of the points $U$ and $V$ in step (1d i) the probability that step (1d ii) succeeds is at least
\begin{align*}
	\left( 1 - \frac{1}{\ell^2}\right)^2.
\end{align*}
\end{Lem}

At the end of step (1d) we have two random $\ell$-torsion points of $E$ namely, $P$ and $Q$. The probability
that $Q$ belongs to the cyclic group generated by $P$ is $\frac{\ell}{\ell^2} = \frac{1}{\ell}$. Thus with
high probability we will find in step (1e) two generators for $E[\ell]$.

\begin{Lem}\label{lem_run_step_1:this}
The expected running time of Step 1 is $O(\ell^{4+o(1)} \log^{2+o(1)} q)$ and
$O(\ell^6 \log^3 q)$ if fast multiplication is not used.
\end{Lem}
\begin{Proof}
The finite field $\mathbb{F}_{q^n}$ can be constructed by picking an irreducible polynomial of degree $n$.
A randomized method that requires on average $O\bigl((n^2 \log n + n \log q) \log n \log \log n\bigr)$ operations over $\mathbb{F}_q$
is given in \cite{sho94}. Thus the field can be constructed in $O(\ell^{4 + o(1)} \log^{2+o(1)} q)$ time since $n \leq \ell^2$.
Step (1d) requires picking a random point on $E$. We can do this by picking a random
element in $\mathbb{F}_{q^n}$ treating it as the $x$-coordinate of a point on $E$ and solving the resulting quadratic
equation for the $y$-coordinate. Choosing a random element in $\mathbb{F}_{q^n}$ can be done in $O(\ell^2 \log q)$ time.
Solving the quadratic equation can be done probabilistically in $O(\ell^2 \log q)$ field operations. Thus to pick a point on $E$
can be done in $O(\ell^{4+o(1)} \log^{2+o(1)} q)$ time. The computation in steps (1d i -- iv) computes integer multiples of a point on
the elliptic curve, where the integer is at most $q^n$,
and this can be done in $O(\ell^{4+o(1)} \log^{2+o(1)} q)$ time using the repeated squaring method and fast multiplication. 
Shanks's Baby-steps-giant-steps
algorithm for a cyclic group $G$ requires $O(\sqrt{|G|})$ group operations. Thus step (1e) runs in time
$O(\ell^{\frac{5}{2} + o(1)} \log^{1+o(1)} q)$, since the group is cyclic of order $\ell$.
\end{Proof}

Let $C$ be a subgroup
of $E$, V{\'e}lu in \cite{vel71} gives explicit formulas for determining the equation of the isogeny $E \rightarrow E/C$
and the Weierstrass equation of the curve $E/C$. 
We give the formulas when $\ell$ is an odd prime. 
Let $E$ is given by the equation
\begin{align*}
y^2 + a_1 xy + a_3 y = x^3 + a_2 x^2 + a_4 x + a_6.
\end{align*}
We define the following two functions in $\mathbb{F}_q(E)$
for $Q = (x,y)$ a point on $E - \{\mathcal{O}\}$ define
\begin{align*}
g^x(Q) &= 3x^2 + 2a_2x + a_4 - a_1 y \\
g^y(Q) &= -2y - a_1 x - a_3,
\end{align*}
and set
\begin{align*}
t(Q) &= 2g^x(Q) - a_1 g^y(Q)\\
u(Q) &= (g^y(Q))^2\\
t	&= \sum_{Q \in (C - \{\mathcal{O}\})} t(Q) \\
w	&= \sum_{Q \in (C - \{\mathcal{O}\})} (u(Q) + x(Q)t(Q)).\\
\end{align*}
Then the curve $E/C$ is given by the equation 
\begin{align*}
Y^2 + A_1XY + A_3Y = X^3 + A_2X^2 + A_4X + A_6
\end{align*}
where
\begin{align*}
A_1 &= a_1, A_2 = a_2, A_3 = a_3, \\
A_4 &= a_4 - 5t, A_6 = a_6 - (a_1^2 + 4a_2)t - 7w.
\end{align*}

From the Weierstrass equation of $E/C$ we can easily determine the $\mathbf{j}$-invariant of $E/C$. It is clear that this procedure
can be done using $O(\ell)$ elliptic curve operations for each of the groups $G_i$, $1 \leq i \leq \ell + 1$. Thus
step 2 can be done in $O(\ell^{4+o(1)} \log^{1+o(1)} q )$ time steps. Step 3 requires only $O(\ell)$ field operations and
so the running time of the algorithm is dominated by the running time of steps 1 and 2. Note that
the polynomial obtained at the end of Step 3 $\phi_{\ell}(x,j)$ has coefficients in $\mathbb{F}_{p^2}[x]$
since all the curves $\ell$-isogenous to $E$ are supersingular and hence their $\mathbf{j}$-invariants belong to $\mathbb{F}_{p^2}$. 
In summary, we have the following:

\begin{Thm} \label{thm_local:this} Algorithm $1$ computes $\phi_{\ell}(x,j) \in \mathbb{F}_{p^2}[x]$, in fact, the
list of roots of $\phi_{\ell}(x,j)$, and has an expected running
time of $O(\ell^{4+o(1)}\log^{2+o(1)} q)$ and $O(\ell^6 \log^3 q)$ without fast multiplication.
\end{Thm}

For our application of Algorithm 1 we will need the dependence of the running time in terms of the quantity $n$. We make
the dependence explicit in the next theorem.

\begin{Thm} \label{thm_local_explicit:this} With notation as above, Algorithm $1$ computes $\phi_{\ell}(x,j) \in \mathbb{F}_{p^2}[x]$ together with
the list of its roots and has an expected running time of $O(n^{2+o(1)}\log^{2+o(1)}q + \sqrt{\ell}n^{1+o(1)}\log^{1+o(1)}q
+ \ell^2 n^{1+o(1)}\log q)$. If fast multiplication is not used then Algorithm $1$ has an expected running time of
$O(n^3 \log^3 q + \sqrt{\ell} n^2 \log^2 q + \ell^2 n^2 \log^2 q)$.
\end{Thm}

In the case of ordinary elliptic curve, step (1) of Algorithm 1 can still be used, once the number of points
on $E/\mathbb{F}_q$ has been determined, by Lemma \ref{lem_deg_bd:this} the degree of the extension, $n$,
is still $O(\ell^2)$. This leads to the following two results:

\begin{Cor} \label{cor_random_tor:this}
If $E/\mathbb{F}_q$ is an elliptic curve, we can pick a random $\ell$-torsion point on $E(\overline{\mathbb{F}}_q)$
in time $O(\ell^{4+o(1)} \log^{2+o(1)}q + \log^{6+o(1)} q)$ and $O(\ell^6 \log^3 q + \log^8 q)$ without fast multiplication.
\end{Cor}

\begin{Cor} \label{cor_random_isog:this}
If $E/\mathbb{F}_q$ is an elliptic curve, we can construct a random $\ell$-isogenous curve to $E$ in
time $O(\ell^{4+o(1)} \log^{2+o(1)}q + \log^{6+o(1)}q)$ and $O(\ell^6 \log^3 q + \log^8 q)$ without fast multiplication.
\end{Cor}

Factoring a degree $d$ polynomial over a finite field $\mathbb{F}_{q^n}$ requires $O(d^2 (n \log q))$ operations over $\mathbb{F}_{q^n}$.
To factor the $\ell$-division polynomial we need an extension of degree roughly $\ell^2$. 
Thus, if we were to factor the $\ell$-division polynomial to generate the isogeny, we would need $O(\ell^{6} \log q)$
operations over a field of degree $\ell^2$ over $\mathbb{F}_q$ which translates to $O(\ell^{8+o(1)}\log^{2+o(1)} q)$ bit operations {\em even}
with fast multiplication.

\section{Computing $\phi_{\ell}(x,y) \mod p$}

In characteristic $p > 2$ there are exactly 
\begin{align*}
S(p) = \biggl\lfloor \frac{p}{12}\biggr\rfloor + \epsilon_p
\end{align*}
supersingular $\mathbf{j}$-invariants where 
\begin{align*}
	\epsilon_p = 0,1,1,2 \text{ if } p \equiv 1,5,7,11 \mod 12.
\end{align*}

In this section we provide an algorithm for computing $\phi_{\ell}(x,y) \mod p$ provided $S(p) \geq \ell+1$.
The description of the algorithm follows:\\

{\bf Algorithm 2}\\
{\bf Input: } Two distinct primes $\ell$ and $p$ with $S(p) \geq \ell + 1$.\\
{\bf Output: } The polynomial $\phi_{\ell}(x,y) \in \mathbb{F}_p[x,y]$.\\
\begin{enumerate}
\item Find the smallest (in absolute value) discriminant $D < 0$ such that $\bigl( \frac{D}{p} \bigr) = -1$.
\item Compute the Hilbert Class polynomial $H_D(x) \mod p$.
\item Let $j_0$ be a root of $H_D(x)$ in $\mathbb{F}_{p^2}$.
\item Set $i = 0$.
\item Compute $\phi_i = \phi_{\ell}(x,j_i) \in \mathbb{F}_{p^2}$ using Algorithm 1.
\item Let $j_{i+1}$ be a root of $\phi_k$ for $k \leq i$ which is not one of $j_0,\cdots, j_i$.
\item If $i < \ell$ then set $i = i+1$ and repeat Step 5.
\item Writing $\phi_{\ell}(x,y) = x^{\ell+1} + \sum_{0 \leq k \leq \ell} p_k(y) x^k$, we have $\ell+1$ systems of equations
			  of the form $p_k(j_i) = v_{ki}$ for $0 \leq k,i \leq \ell$. 
		Solve these equations for each $p_k(y)$, $0 \leq k \leq \ell$.
\item Output $\phi_{\ell}(x,y) \in \mathbb{F}_p[x,y]$.\\
\end{enumerate}

We argue that the above algorithm is correct and analyze the running time. For step $1$, we note
that if $p \equiv 3 \mod 4$, then $D = -4$ works. Otherwise, $-1$ is a quadratic residue and writing (without loss of generality)
$D$ as $-4d$, we are looking for the smallest $d$ such that $\bigl( \frac{d}{p} \bigr) = -1$. A theorem of Burgess (\cite{bur62}) tells us
that $d \ll p^{\frac{1}{4\sqrt{e}}}$, and under the assumption of GRH the estimate of Ankeny (\cite{ank52}) gives $d \ll \log^2 p$.
Computing $H_D(x) \mod p$ can be done in $O(d^2 (\log d)^2)$ time \cite{ll90} \S5.10. Thus step 2 requires $O(\sqrt{p}\log^2 p)$ time, 
and under the assumption of GRH requires $O\bigl(\log^4 p (\log \log p)^2\bigr)$ time.
Since $\bigl( \frac{D}{p} \bigr) = -1$ all the roots of $H_D(x)$ are supersingular $\mathbf{j}$-invariants in characteristic $p$.
$H_D(x)$ is a polynomial of
degree $h(\sqrt{-D})$, the class number of the order of discriminant $D$, and this is $\ll |D|^{\frac{1}{2}+\epsilon}$. 
Finding a root of $H_D(x) \in \mathbb{F}_{p^2}$ can be done in $O(d^{1+\epsilon} \log^{2+o(1)} p)$ time using probabilistic factoring
algorithms, where $d = |D|$. 
The graph with supersingular $\mathbf{j}$-invariants over charactertistic $p$ as vertices and $\ell$-isogenies as edges is connected (see \cite{mes86}),
consequently, we will always find a $\mathbf{j}$-invariant in step 6 that is not one of $j_0,\cdots, j_{i}$. Thus the loop in steps $(5)\cdots(7)$
is executed exactly $\ell+1$ times under the assumption that $S(p) \geq \ell+1$. Even though Algorithm $1$ requires $\tilde{O}(\ell^4 \log^2 q)$ 
time\footnote{We use the soft-Oh $\tilde{O}$ notation when we ignore factors of the form $\log \ell$ or $\log\log p$.} 
in the worst case, we will argue that, in fact, 
for {\em all} of the iterations of the loop it actually runs in $\tilde{O}(\ell^3 \log^2 q)$ time.

\begin{Lem} Let $E$ be a supersingular elliptic curve defined over $\mathbb{F}_{p^2}$.
Then the extension degree $$[\mathbb{F}_{p^2}(E[\ell]) : \mathbb{F}_{p^2}]$$
divides $6(\ell-1)$.
\end{Lem}
\begin{Proof}
Let $E/\mathbb{F}_{p^2}$ be a supersingular curve and let $t$ be the trace of Frobenius. Then the Frobenius map $\phi$
satisfies
\begin{align*}
	\phi^2 - t\phi + p^2 = 0
\end{align*}
with $t = 0$ or $\pm p$ or $\pm 2p$.
Suppose $t = \pm 2p$, then the Frobenius acts as multiplication by $\pm p$ on the curve $E$. 
Thus $\phi^{\ell-1}$ acts trivially on $E[\ell]$, and the $\ell$-torsion points
are defined over an extension of degree dividing $\ell-1$. 
If $t = 0$, then $\phi^2 = -p^2$ and so 
$\phi^{2(\ell-1)}$ acts trivially on the $\ell$-torsion. Thus $E[\ell]$ is defined over an extension
of degree dividing $2(\ell-1)$. If $t = \pm p$, then $\phi^3 = \pm p^3$ and consequently $\phi^{3(\ell-1)}$
acts trivially on the $\ell$-torsion of the curve. Thus the $\ell$-torsion is defined over an extension
of degree dividing $3(\ell-1)$. Thus in all cases the $\ell$-torsion of the curve is defined over
an extension of degree dividing $6(\ell-1)$.
\end{Proof}

Thus, the loop in steps $(5) \cdots (7)$, Algorithm 1 can be run with the quantity $n = 6(\ell-1)$. 
For this value of $n$ Algorithm 1 runs in expected time $O(\ell^{3+o(1)} \log^{2+o(1)} p)$ 
and so the loop runs in expected time $O(\ell^{4+o(1)} \log^{2+o(1)} p)$.\\

Writing the modular polynomial $\phi_{\ell}(x,y)$ as $x^{\ell+1} + \sum_{0 \leq k \leq \ell} p_k(y) x^{k}$, we know that 
$p_0(y)$ is monic of degree $\ell+1$ and $\deg(p_k(y)) \leq \ell$ for $1 \leq k \leq \ell$. Thus at the end of the loop in steps $(5) \cdots (7)$
we have enough information to solve for the $p_k(y)$ in step (8).
We are solving $\ell+1$ systems of equations, each
requiring an inversion of a matrix of size $(\ell+1) \times (\ell+1)$. This can be done in $O(\ell^4\log^{1+o(1)} p)$ time.
Since the polynomial $\phi_{\ell}(x,y) \mod p$ is the reduction of the classical modular polynomial, a polynomial
with integer coefficients, the polynomial that we compute has coefficients in $\mathbb{F}_p$. Thus we have proved 
the following theorem:

\begin{Thm} \label{thm_mod_p_comp:this}
Given $\ell$ and $p$ distinct primes such that $S(p) \geq \ell+1$, Algorithm 2 computes $\phi_{\ell}(x,y) \in \mathbb{F}_p[x,y]$
in expected time $O(\ell^{4+o(1)}\log^{2+o(1)} p + \log^4 p \log\log p)$ under the assumption of GRH and
in $O(\ell^5 \log^3 p)$ time without fast multiplication.
\end{Thm}

Hence, we can compute $\phi_{\ell}(x,y)$ modulo a prime $p$ in $\tilde{O}(\ell^{4} \log^{2} p + \log^4 p)$ time if $p \geq 12 \ell + 13$.
If $p < 12\ell + 13$, we could still use the algorithm with ordinary elliptic curves and this would lead
to a running time with the dependence on $\ell$ being $\ell^5$. Furthermore, we would not need the GRH since
it was needed only to determine a supersingular curve in characteristic $p$. However, this is not very efficient.

\begin{Rem} If one is allowed to pick the prime $p$, such as would be the case if we are computing $\phi_{\ell}$ over
the integers using the Chinese Remainder Theorem combined with this method, then one can eliminate the assumption of GRH in the above theorem. For example,
for primes $p \equiv 3 \mod 4$ the ${\bf j}$-invariant $1728$ is supersingular. Thus in step (3) of Algorithm 2, we can
start with $j_0 = 1728$ for any such prime. Hence we do not need the GRH to bound $D$ in the analysis of the running time
of the algorithm.
\end{Rem}

{\bf Acknowledgements:} We would like to thank Igor Shparlinski and Steven Galbraith for useful suggestions and
comments on an earlier draft of the paper. We would like to
thank the anonymous referee for suggestions that improved the presentation 
of the paper. The first author would like to thank Microsoft Research for providing a stimulating
research environment for this work.

\section*{Appendix}
Elkies in \cite{elk98} (see \S3 of that paper) claims that the usual method of computing the $\ell$-th modular polynomial runs in time 
$O(\ell^{4+\epsilon})$. However, there is a subtle error in the analysis. We argue that in fact, the running
time of this algorithm is $\Theta(\ell^{\frac{9}{2} + \epsilon})$. The first stage of the algorithm
invovles computing the first $\ell^2+O(\ell)$ Fourier coefficients of the powers of the $\mathbf{j}$-function, namely
$\mathbf{j}$, $\cdots$, $\mathbf{j}^{\ell}$. This (as Elkies points out) can be done in $O(\ell^{3+\epsilon})$ arithmetic operations.
The problem comes when we study the running time in terms of {\em bit-operations}. To analyze this we need to
study the bit-sizes of the numbers that are handled by the algorithm. 
While it is true that the $n$-th Fourier
coefficient of $\mathbf{j}$ grows as $e^{O(\sqrt{n})}$, we need to also compute 
the Fourier coefficients of powers of $\mathbf{j}$
and they grow faster (essentially because they have a higher order pole at $\infty$). 
In \cite{mah74} 
an upper bound of the form $\exp(4\pi(\sqrt{(n+k)k})$ is proven for the $n$-th Fourier coefficient of $\mathbf{j}^k$. We show
that the upper bound is quite close to the true magnitude of the Fourier coefficients below.
Let $c(n)$ denote the $n$-th Fourier coefficient of the $\mathbf{j}$-function. It is well known that $c(n)$ are
all positive integers and that (see \cite{pet32})
\begin{align*}
c(n) \sim \frac{e^{4\pi \sqrt{n}}}{\sqrt{2} n^{3/4}}.
\end{align*}

The $n$-th Fourier coefficient of $\mathbf{j}^k$ is given by
\begin{align*}
	\sum_{a_1 + a_2 + \cdots + a_k = n} c(a_1) c(a_2) \cdots c(a_k).
\end{align*}
Clearly, there is at least one partition of $n$ (into $k$ parts) 
where each of the parts, $a_i$, are $\geq \frac{n}{ck}$, where $c > 1$ is a constant.
The asymptotic formula for $c(n)$ gives us that
\begin{align*}
\sum_{a_1 + a_2 + \cdots + a_k = n} c(a_1) c(a_2) \cdots c(a_k) \gg \left(e^{\sqrt{n/ck} - O(\log n)}\right)^k = e^{\Omega(\sqrt{nk})},
\end{align*}
as long as $k$ and $n$ vary such that the ratio $n/k$ goes to infinity.
Thus a lower bound for the rate of growth of the $n$-th Fourier coefficient
of $\mathbf{j}^k(z)$ is $e^{\Omega(\sqrt{kn})}$. 
Hence to compute the first $\ell^2 + O(\ell)$ Fourier coefficients of the powers 
$\mathbf{j}^k$, for $1 \leq k \leq \ell$, the bit length of the numbers
involved is $O(\ell\sqrt{\ell})$ rather than the estimate $O(\ell \log \ell)$ used in \cite{elk98}. 
Thus this stage of the algorithm already requires $\Theta(\ell^{\frac{9}{2}+\epsilon})$
time using fast multiplication. 
The dependence of the running time on $\ell$ for our method, on the other hand, is $O(\ell^{4+\epsilon})$ if fast multiplication is
used.
The error in \cite{elk98} stems from the fact that a bound on the size of the coefficients
of $\phi_{\ell}(x,y)$ was used to bound the bit sizes. These coefficients are somewhat smaller, their absolute value
being bounded by $e^{O(\ell \log \ell)}$ (see \cite{coh84}).
\end{document}

%% file: MTemplate.tex
\usepackage{fullpage}
\usepackage{amssymb}
\usepackage{amsmath}
\usepackage{amsfonts}
\usepackage{amsxtra}
\usepackage{amscd}
\usepackage{amsthm}
\usepackage{epsfig}
\usepackage{enumerate}
\usepackage{float}
\usepackage{xy}

\xyoption{all}

\theoremstyle{plain} \newtheorem{Thm}{Theorem}[section]
\theoremstyle{plain} \newtheorem{Cor}[Thm]{Corollary}
\theoremstyle{plain} \newtheorem{Lem}[Thm]{Lemma}
\theoremstyle{plain} 
\theoremstyle{definition} 
\theoremstyle{definition} 
\theoremstyle{remark} \newtheorem{Rem}[Thm]{Remark}
\newenvironment{Proof}{{\bf Proof :}}{$\Box$\\}

\newcommand{\Gal}{\mbox{$\mathrm{Gal}$}}

\newcommand{\GL}{\mbox{$\mathrm{GL}$}}

\newcommand{\newfontobj}[2]{
  \newcommand{#1}[1]{
    \expandafter\def\csname##1\endcsname{{#2 ##1}}}}

\newfontobj{\class}{\rm}
\newfontobj{\lang}{\bf}
\newfontobj{\oper}{\rm}

\class{PSPACE}		
\class{AM}		
\class{MA}
\class{NP}
\class{UP}
\class{P}
\class{RP}
\class{BPP}
\class{DTIME}
\class{ZPP}
\class{EXPSPACE}
\class{coNP}
\class{coRP}
\class{coAM}
\class{PH}
\class{co}
\class{BP}	

\lang{HN} 		
\lang{SAT}
\lang{PIM}
\lang{CSG}
\lang{ESC}
\lang{IP}
\lang{PRIMES}

\oper{Pr}    
\oper{E}     
\oper{Ord}   
\oper{li}

\floatstyle{ruled}
\newfloat{Algorithm}{thp}{loa}[section]


%% file: Modular_Poly_Final_LMSJCM.bbl
\begin{thebibliography}{cl04}
\bibitem[ALV03]{alv03}
Agashe, A.; Lauter, K.; Venkatesan, R.;
{\em Constructing Elliptic Curves with a known number of points over a prime field}, in
Lectures in honour of the 60th birthday of Hugh Cowie Williams, Fields Institute Communications Series, {\bf 42}, 1-17, 2003.

\bibitem[Ank52]{ank52}
Ankeny, N., C.; {\em The least quadratic non-residue}, Annals of Math., (2), {\bf 55}, 65-72, 1952.

\bibitem[BSS99]{bss99}
Blake, I.; Seroussi, G.; Smart, N.;
{\em Elliptic Curves in Cryptography}, Lond. Math. Soc., Lecture Note Series, {\bf 265}, Cambridge University
Press, 1999.


\bibitem[Bur62]{bur62}
Burgess, D., A.; {\em On character sums and primitive roots}, Proc. London Math. Soc., {\bf III}, {\bf 12}, 179-192, 1962.

\bibitem[CNST98]{cnst98}
Chao, J.; Nakamura, O.; Sobataka, K.; Tsujii, S.;
{\em Construction of secure elliptic cryptosystems using CM tests and liftings},
Advances in Cryptology, ASIACRYPT'98 (Beijing), Lecture Notes in Computer Science, 1514,
Springer-Verlag, Berlin, 1998.

\bibitem[Coh84]{coh84}
Cohen, P.; 
{\em On the coefficients of the transformation polynomials for the elliptic modular function},
Math. Proc. Cambridge Philos. Soc., {\bf 95}, 389-402, 1984.

\bibitem[Elk98]{elk98}
Elkies, Noam; {\em Elliptic and modular curves over finite fields and related computational issues},
in Computational Perspectives on Number Theory: Proceedings of a Conference in Honor of A.O.L. Atkin 
(D.A. Buell and J.T. Teitelbaum, eds.), AMS/International Press, 21-76, 1998. 


\bibitem[LL90]{ll90}
Lenstra, A., K.; Lenstra, H., W., Jr.;
{\em Algorithms in Number Theory}, Handbook of Theoretical Computer Science, Vol. {\bf A}. Elsevier, 673-715, 1990.

\bibitem[Mah74]{mah74}
Mahler, K.;
{\em On the coefficients of transformation polynomials for the modular function}.
Bull. Austral. Math. Soc., {\bf 10}, 197-218, 1974.

\bibitem[Mes86]{mes86}
Mestre, J.-F.;
{\em La m{\'e}thode des graphes. Exemples et applications},
Proceedings of the international conference on class numbers and fundamental units of algebraic number fields,
Nagoya Univ., Nagoya, 217-242, 1986.

\bibitem[Pet32]{pet32}
Petersson, H.; {\em {\"U}ber die Entwicklungskoeffizienten der automorphen formen}, 
Acta Math., {\bf 58}, 169-215, 1932.

\bibitem[Sch87]{sch87}
Schoof, R.; {\em Nonsingular Plane Cubic Curves over Finite Fields},
J. Combinatorial Theory, {\bf 46}, 2, 183-208, 1987.

\bibitem[Sho94]{sho94}
Shoup, V.; {\em Fast construction of irreducible polynomials over finite fields}, J. Symbolic Computation, {\bf 17}, 371-391, 1994.

\bibitem[Sil86]{sil86} 
Silverman, J. H.; {\em The Arithmetic of Elliptic Curves}, Graduate Texts in Mathematics, {\bf 106}, 
Springer-Verlag, 1986.

\bibitem[Vel71]{vel71}
V{\'e}lu, J.; {\em Isog{\'e}nies entre courbes elliptiques}, C. R. Acad. Sc. Paris, {\bf 273}, 238-241, 1971.

\end{thebibliography}
